\documentclass[11pt,amstex,amssymb,reqno]{amsart}
\usepackage{amsmath,amsthm,amsfonts,amssymb,amscd}
\usepackage[english]{babel}
\usepackage[utf8]{inputenc}
\usepackage{float}
\usepackage[all]{xy}
\usepackage{float}
\usepackage{listings}
\usepackage{subfig}
\usepackage{listings}
\usepackage{xcolor}
\definecolor{codegreen}{rgb}{0,0.6,0}
\definecolor{codegray}{rgb}{0.5,0.5,0.5}
\definecolor{codepurple}{rgb}{0.58,0,0.82}
\definecolor{backcolour}{rgb}{0.95,0.95,0.92}

\lstdefinestyle{mystyle}{
    backgroundcolor=\color{backcolour},
    commentstyle=\color{codegreen},
    keywordstyle=\color{magenta},
    numberstyle=\tiny\color{codegray},
    stringstyle=\color{codepurple},
    basicstyle=\ttfamily\footnotesize,
    breakatwhitespace=false,
    breaklines=true,
    captionpos=b,
    keepspaces=true,
    numbers=left,
    numbersep=5pt,
    showspaces=false,
    showstringspaces=false,
    showtabs=false,
    tabsize=2
}

\usepackage{amssymb,amsthm,amsmath,psfrag,mathrsfs}
\usepackage[dvips]{graphicx}
\usepackage{xcolor}
\usepackage{comment}

\newtheorem{algorit}{Algoritmo}[section]

\makeatletter
\def\author@andify{%
  \nxandlist {\unskip ,\penalty-1 \space\ignorespaces}%
    {\unskip {} \@@and~}%
    {\unskip \penalty-2 \space \@@and~}%
}
\makeatother

\title[Numerical approximation of the insitu combustion model using the nonlinear mixed complementarity method]{Numerical approximation of the insitu combustion model using the nonlinear mixed complementarity method}

\author{Julio C\'esar Agust\'in Sangay}
\address{J. Sangay. Universidad Privada del Norte, Trujillo-Per\'u}
\email{julio.agustin@upn.edu.pe}

\author{Alexis Rodriguez Carranza}
\address{A. Rodriguez. Universidad Privada del Norte, Trujillo-Per\'u}
\email{alexis.rodriguez@upn.edu.pe}

\author{George J. Bautista}
\address{Faculty of Engineering, Civil Engineering Professional School, Technological University of the Andes(UTEA), Sede Abancay, Av. Per\'u 700, Apur\'imac, Peru.}
\email{gbautistas@utea.edu.pe}

\author{Juan Carlos Ponte Bejarano}
\address{Universidad Tecnol\'ogica del Peru, Trujillo-Per\'u}
\email{C25226@utp.edu.pe}
\author{Jos\'e Luis Ponte Bejarano}
\address{J. Ponte Universidad Tecnol\'ogica del Peru, Trujillo-Per\'u}
\email{c24064@utp.edu.pe}
\author{Eddy Cristiam Miranda Ramos}
\address{E. Cristiam Universidad Nacional de Trujillo}
\email{ emirandar@unitru.edu.pe}

\begin{document}

\begin{abstract}
In this work, we will study a numerical method that allows finding an approximation of the exact solution for a in-situ combustion model using the nonlinear mixed complementary method, which is a variation of the Newton's method for solving nonlinear systems based on an implicit finite difference scheme and a nonlinear algorithm mixed complementarity, FDA-MNCP. The method has the advantage of provide a global convergence in relation to the finite difference method and method of Newton that only has local convergence. The theory is applied to model in-situ combustion, which can be rewritten in the form of mixed complementarity also we do a comparison with the FDA-NCP method.
\end{abstract}
\maketitle

\tableofcontents

\section{Introduction}
several mathematical models in different disciplines such as engineering, physics, economics and other sciences study partial differential equations of the parabolic type. These models can lead to the problem of mixed complementarity, that is, the case 
several mathematical models in different disciplines, such as engineering, physics, economics and other sciences study partial differential equations of the parabolic type. These models can lead to the problem of mixed complementarity, that is, the case of the in-situ combustion model, which will be our model in this work. Other applications of complementarity problems are described in \cite{ferris}.

Since the goal is to find an approximation of the analytical solution, we will develop a numerical method that will allow us to achieve our goal. This technique will be applied to the simple in-situ combustion model which will be written as a mixed complementarity problem. 

The in-situ combustion model is a particular case of the model treated in \cite{chapiro}.  
In this case, the model considers the injection of air into a porous medium containing solid fuel and consists of a system of two nonlinear parabolic differential equations. 

The contribution of the work is the study of the simple in-situ combustion model and simulations for the proposed model applying the Crank-Nicolson method and the FDA-MNCP algorithm \cite{sandro}. 

We will present the results that indicate that the sequence of feasible points generated is contained in a feasible region and we will verify that the directions obtained are feasible and descending for a function associated with the complementarity problem and we will also see the proof of global convergence for the FDA-MNCP following the feat for FDA-NCP\cite{sandrofdancp}. 
We apply the FDA-MNCP method to the in-situ combustion problem, describe the discretization procedure using the finite difference technique for the mixed complementarity problem associated with the problem, and also present the numerical results and the corresponding error analysis with the comparison with the FDA-NCP method. Finally, we will present some conclusions.

\section{Physical problem modeling}
The model studies one-dimensional flows, with a combustion wave in the case when the oxidant (air with oxygen) is injected into a porous medium. Initially the medium contains a fuel that is essentially immobile, does not vaporize and the amount of oxygen is unlimited. This is the case for solid or liquid fuels with low saturations.
As in \cite{angel}, we study the simplified model where:

\begin{itemize}
\item[ $\bullet$ ] A small part of the available space is occupied by fuel.
\item[ $\bullet$ ] Porosity changes in the reaction are negligible.
\item[ $\bullet$ ] Temperature of solid and gas are the same (local thermal equilibrium).
\item[ $\bullet$ ] Gas velocity is constant.
\item[ $\bullet$ ] Heat loss is negligible.
\item[ $\bullet$ ] Pressure variations are small compared to the prevailing pressure.
\end{itemize}
The model has temporal $t$ and spatial $x$ coordinates that include the heat balance equation, the molar balance equation for immobile fuels, and the ideal gas law.
\begin{eqnarray}{\label{S1}}
C_{m} \frac{\partial T}{\partial t} + \frac{\partial}{\partial x}(C_{g} \rho u(T - T_{res})) &=& \lambda \frac{\partial^{2} T}{\partial x^{2}} + Q_{r} W_{r}, \\
\frac{\partial \rho_{f}}{\partial t} &=& - \mu_{f} W_{r}, \\
\rho &=& \frac{P}{TR},
\end{eqnarray}
where $T [K]$ is the temperature, $\rho [\frac{mol}{m^{3}}]$ is the molar density of the gas, and $\rho_{f} [\frac{mol}{m^{3}}]$ is the molar concentration of the immobile fuel. The set of parameters along with their typical values are given in the Table\ref{tabeladecombustaoinsitu}.

\begin{table}[hbt]
\begin{center}
\begin{tabular}{|p{1.4cm}|p{8.6cm}|p{1.2cm}|p{2.0cm}|}   \hline  
\textit{Symbol} & \textit{Physical Quantity} & \textit{Value} & \textit{Unit} \\ \hline
$T_{res}$        & Initial reservoir temperature & 273             & $[K]$   \\ \hline
$C_{m}$          & Heat capacity of porous medium & $2 \cdot 10^{6}$     & $[J/m^{3}K]$    \\ \hline
$c_{g}$          & Heat capacity of gas &  27.42          & $[J/mol K]$   \\ \hline
$\lambda$        & Thermal conductivity of porous medium & 0.87            & $[J/(m s K)]$   \\ \hline
$Q_{r}$          & Enthalpy of the still fuel in $T_{res}$          & $4\cdot 10^{5} $       & $[J/mol]$   \\ \hline
$u_{inj}$        & Darcy speed for gas injection $(200 m/dia)$    & $0.0023$        & $[m/s]$   \\ \hline
$E_{r}$          & Activation energy & $58000$          & $[J/mol]$   \\ \hline
$K_{p}$          & Pre-exponential parameter & 500             & $1/s$    \\ \hline
$R$              & Ideal gas constant & 8.314           & $[J/mol K]$    \\ \hline
$P$              & Prevailing pressure $(1 atm)$                       & 101325          & $[Pa]$    \\ \hline
$\rho^{res}_{f}$ & Initial molar density of fuel & 372             & $[mol/m^{3}]$    \\ \hline
\end{tabular}
\end{center}
\caption{Dimensional parameters for in-situ combustion and their typical values \cite{angel}.} 
\label{tabeladecombustaoinsitu}
\end{table} 
As in \cite{angel}, if we consider for simplicity $\mu_{f}=\mu_{g}=\mu_{0}=1$ and the amount of oxygen is unlimited, the reaction ratio $W_{r}$ is taken as:
\begin{equation}
W_{r} = k_{p} \rho_{f} \exp \left(\frac{-E_{r}}{RT}\right).
\end{equation}
The variables to be found are the temperature $(T)$ and the molar concentration of the immobile fuel $(\rho_{f})$. Since the equations are not dimensioned, we do as \cite{angel}, to obtain the dimensionless form:
\begin{eqnarray} \label{5.5}
\frac{\partial \theta}{\partial t} + u \frac{\partial (\rho \theta)}{\partial x} &=& \frac{1}{P_{e_{T}}} \frac{\partial^{2} \theta}{\partial x^{2}} + \Phi (\theta, \eta). \\
\frac{\partial \eta}{\partial t} &=& \Phi (\theta, \eta).
\end{eqnarray}
$$
\mbox{where:} \ \ \ \  \rho = \frac{\theta_{0}}{\theta + \theta_{0}}, \ \ \ \ \ \ 
\Phi = \beta (1- \eta) \exp\left(\frac{-E}{\theta + \theta_{0}}\right).
$$
with the dimensionless constants:
\begin{equation} \label{constantesadimensionais}
P_{E_{T}} = \frac{x^{\star}}{\lambda \Delta T^{\star}}, \ \ \beta = \rho_{f}^{\star} k_{p} Q_{r}, \ \ E = \frac{E_{r}}{R \Delta T^{\star}}, 
\ \ \theta_{0} = \frac{T_{res}}{\Delta T^{\star}}, \ \ u = \frac{u_{inj} t^{\star}}{x^{\star}}.
\end{equation}

Here $P_{e_T}$ is the Peclet number for thermal diffusion, $u$ becomes the dimensionless thermal wave velocity, $E$ is a escalated activation energy and $\theta_0$ is a scaled reservoir temperature. With a reservoir initial conditions:

$$
t=0;\quad x\geq 0:\quad \theta=0,\quad \eta=0
$$
and injections condition:
$$
t\geq 0;\quad x= 0:\quad \theta=0,\quad \eta=1.
$$

\section{Description of the FDA-MNCP method for the simple in-situ combustion model}
We now describe in detail the finite difference scheme for the in-situ combustion model for the FDA-MNCP method. For this we use the mesh applying the Crank-Nicolson method to approximate the spatial derivatives in each case, i.e.:

\begin{eqnarray} \label{malha-insitu}
\partial_{t} \theta(x_{m},t_{n+\frac{1}{2}}) & = & \frac{\theta^{n+1}_{m}-\theta^{n}_{m}}{k}. \label{diferencia progresiva, adelantada o posterior-insitu} \\
\partial_{xx} \theta(x_{m},t_{n+\frac{1}{2}}) & = & \frac{\theta^{n+1}_{m+1}-2 \theta^{n+1}_{m} + \theta^{n+1}_{m-1}}{2h^{2}} + \frac{\theta^{n}_{m+1}-2 \theta^{n}_{m} + \theta^{n}_{m-1}}{2h^{2}}. \label{diferencia central-insitu} \\
\partial_{x} F(\theta(x_{m},t_{n+\frac{1}{2}})) & = & \frac{F^{n+1}_{m+1}-F^{n+1}_{m-1}}{4h} + \frac{F^{n}_{m+1}-F^{n}_{m-1}}{4h}. \\
\Phi(\theta(x_{m},t_{n+\frac{1}{2}})) & = & \frac{\Phi^{n+1}_{m} + \Phi^{n}_{m}}{2}. \label{phi-insitu}
\end{eqnarray}

Considering the Dirichlet conditions at the point $x_0$:

$$
\theta(x_0,t)=0,\quad\quad \eta(x_0,t)=1,
$$

and Neumann conditions at the point $x_M$
$$
\frac{\partial\theta}{\partial x}(x_M,t)=0,\quad\quad \frac{\partial\eta}{\partial x}(x_M,t)=0,
$$
given in\cite{31 angel}, therefore, the value at $x_0$ is known at all times but not at $x_M$, and
\begin{equation}{ \label{direchelt-insitu}}
\theta^{n+1}_{0}=\theta^{n}_{0}, \quad \eta^{n+1}_{0}=\eta^{n}_{0},\quad \text{for all $n\in\mathbb{N}$}
\end{equation}
The boundary condition on $x_M$ gives:
$$
\frac{\partial\theta}{\partial x}(x_M,t)=0\Longrightarrow \frac{\theta^{n}_{M+1}-\theta^{n}_{M-1}}{2h}=0
$$
therefore,
\begin{equation}{ \label{5.18ang-insitu}}
\theta^{n}_{M+1}=\theta^{n}_{M-1}\Longrightarrow F^{n}_{M+1}=F^{n}_{M-1},\quad \text{for all $n\in\mathbb{N}$}
\end{equation}
given shaping the FDA-MNCP method:

\begin{eqnarray} \label{F en situ}
\theta \geq 0 ; \ \frac{\partial \theta}{\partial t} + u \frac{\partial (\rho \theta)}{\partial x} - \frac{1}{P_{e_{T}}} \frac{\partial^{2} \theta}{\partial x^{2}} - \Phi (\theta, \eta) \geq 0. \\
\frac{\partial \eta}{\partial t}- \Phi(\theta, \eta)=0. \label{Q en situ} \end{eqnarray}
To obtain the discrete relations of(\ref{F en situ}), we replace (\ref{malha-insitu}-\ref{phi-insitu}) in(\ref{F en situ}) to obtain:

\begin{eqnarray}\label{discrete-F ensitu}
-2 \mu H \theta^{n+1}_{m-1}+(4+4 \mu H) \theta^{n+1}_{m}-2 \mu H \theta^{n+1}_{m+1} + \lambda [F^{n+1}_{m+1}-F^{n+1}_{m-1}] - 2k \Phi^{n+1}_{m} \ \ \geq \\ \nonumber
\ \ \  \ \ \ \ \ \ \ \ \geq \ \ 2 \mu H \theta^{n}_{m-1}+(4-4 \mu H) \theta^{n}_{m}+2 \mu H \theta^{n}_{m+1} + \lambda [F^{n}_{m+1}-F^{n}_{m-1}] + 2k \Phi^{n}_{m},
\end{eqnarray}
where $\lambda=\frac{h}{h}$ and $\mu=\frac{k}{h^^2}$. The scheme is valid for $m=1,2,\dots,M$ of the points whose values are not known.
At the boundary points we have that for $m = 1$, substituting (\ref{direchelt-insitu}) in (\ref{discrete-F ensitu}) we obtain:
\begin{eqnarray}\label{F com m=1 ensitu}
(4+4 \mu H) \theta^{n+1}_{1}-2 \mu H \theta^{n+1}_{2} + \lambda [F^{n+1}_{2}-F^{n+1}_{0}] - 2k \Phi^{n+1}_{1} \ \ \geq \\ \nonumber
\ \ \geq \ \ (4-4 \mu H) \theta^{n}_{1}+2 \mu H \theta^{n}_{2} - \lambda [F^{n}_{2}-F^{n}_{0}] + 2k \Phi^{n}_{1} + 4 \mu H \theta^{n}_{0},
\end{eqnarray}
for $n=M$, we replace (\ref{5.18ang-insitu}) in (\ref{discrete-F ensitu}) for obtain:
\begin{eqnarray}\label{F com m=M ensitu}
-4 \mu H \theta^{n+1}_{M-1} + (4+4 \mu H) \theta^{n+1}_{M} - 2k \Phi^{n+1}_{M} \geq  4 \mu H \theta^{n}_{M-1} + (4-4 \mu H) \theta^{n}_{M} + 2k \Phi^{n}_{M}.
\end{eqnarray}
Then (\ref{discrete-F ensitu}) is valid for all $m=2,\dots,M-1$ and joining the expressions (\ref{F com m=1 ensitu}) and (\ref{F com m=M ensitu})we obtain the following inequality in the variable $\theta^{n+1}$
\begin{eqnarray}\label{Fdiscre ensitu}
G^{n}(\theta^{n+1}, \eta^{n+1}) =  A \theta^{n+1} + \lambda P(\theta^{n+1}, \eta^{n+1}) - 2k \Phi(\theta^{n+1}, \eta^{n+1})-LD(\theta^{n}, \eta^{n}) \geq 0,
\end{eqnarray}
where $LD = B \theta^{n} - \lambda P(\theta^{n}, \eta^{n}) + 2k \Phi(\theta^{n}, \eta^{n}) + UR$ is known at every instant of time. 

Furthermore:

\begin{equation} \label{A ensitu}
\hspace{-0.26cm} A=\left[
\begin{array}{cccccccc}
\! 4+4 \mu H & -2 \mu H  &    0      &   0      &  \cdots  &     0    &    0     &    0   \\
\! -2 \mu H  & 4+4 \mu H & -2 \mu H  &   0      &  \cdots  &     0    &    0     &    0   \\
\!   0      &  -2 \mu H & 4+4 \mu H & -2 \mu H &          &     0    &    0     &    0   \\
 \vdots   &           &           &          &  \ddots  &          &          &  \vdots \\
\!     0    &    0      &    0      &    0     &          & -2 \mu H &4+4 \mu H & -2 \mu H \\
\!     0    &    0      &    0      &    0     &          &   0      & -4 \mu H & 4+4 \mu H      
     \end{array}
\right] \! ,
\end{equation}

\begin{equation} \label{B ensitu}
B=\left[
 \begin{array}{cccccccc}
4-4 \mu H &  2 \mu H  &    0      &   0      &  \cdots  &     0    &    0     &    0   \\
 2 \mu H  & 4-4 \mu H &  2 \mu H  &   0      &  \cdots  &     0    &    0     &    0   \\
   0      &   2 \mu H & 4-4 \mu H &  2 \mu H &          &     0    &    0     &    0   \\
 \vdots   &           &           &          &  \ddots  &          &          &  \vdots \\
     0    &    0      &    0      &    0     &          &  2 \mu H &4-4 \mu H &  2 \mu H \\
     0    &    0      &    0      &    0     &          &   0      &  4 \mu H & 4-4 \mu H      
     \end{array}
\right],
\end{equation}

\begin{equation}\label{P e phi ensitu}
P^{n}=P(\theta^{n})=\left[
 \begin{array}{ccc}
F^{n}_{2}    -       F^{n}_{0}  \\
F^{n}_{3}    -       F^{n}_{1}  \\
F^{n}_{4}    -       F^{n}_{2}  \\
           \vdots       \\
F^{n}_{M}    -       F^{n}_{M-2}  \\
              0        
     \end{array}
\right], \ \
\Phi^{n}=\Phi(\theta^{n})=\left[
 \begin{array}{ccc}
\Phi^{n}_{1} \\
\Phi^{n}_{2} \\
\Phi^{n}_{3} \\
  \vdots \\
\Phi^{n}_{M-1}  \\ 
\Phi^{n}_{M} \\
     \end{array}
\right], \ \ \
\end{equation}

\begin{equation} 
\theta^{n}=\left[
 \begin{array}{ccc}
\theta^{n}_{1} \\
\theta^{n}_{2} \\
\theta^{n}_{3} \\
  \vdots \\
\theta^{n}_{M-1}  \\ 
\theta^{n}_{M} \\
     \end{array}
\right], \ \ 
UR=\left[
 \begin{array}{ccc}
4 \mu H \theta^{n}_{0} \\
0 \\
0 \\
  \vdots \\
0  \\ 
0 \\
     \end{array}
\right]
\end{equation}
where  $A,B \in \mathbb{R}^{M \times M}; \ \theta^{n}, P^{n}, \Phi^{n} \in \mathbb{R}^{M}.$
Similarly, to obtain the discrete form of (\ref{Q en situ}) we replaced (\ref{malha-insitu}-\ref{phi-insitu}) in  (\ref{Q en situ}) to obtain:

\begin{eqnarray} \label{discrete-Q ensitu}
Diag(2) \eta^{n+1}_{m} - k \Phi^{n+1}_{m} = Diag(2) \eta^{n}_{m} + k \Phi^{n}_{m}.
\end{eqnarray}
The difference scheme is valid for all $m=1,2,\dots,M$ of the points whose values are not known.

At the boundary points we have that for $m=1$, we replace (\ref{direchelt-insitu}) in  (\ref{discrete-Q ensitu}) we obtain:
\begin{equation}\label{Q com m=1 ensitu}
2 \eta^{n+1}_{1} - k \Phi^{n+1}_{1} = 2 \eta^{n}_{1} + k \Phi^{n}_{1},
\end{equation}

for $m=M$, we replaced (\ref{5.18ang-insitu}) in (\ref{discrete-Q ensitu}) to obtain:

\begin{eqnarray}\label{Q com m=M ensitu}
2 \eta^{n+1}_{M} - k \Phi^{n+1}_{M} = 2 \eta^{n}_{M} + k \Phi^{n}_{M}.
\end{eqnarray}

Therefore, (\ref{discrete-Q ensitu}) is valid for all $m=2,\dots,M-1$ and joining the expressions (\ref{Q com m=1 ensitu}) and (\ref{Q com m=M ensitu}) we obtain the following inequality in the variable $\eta^{n+1}$;

\begin{equation} \label{Qdiscre ensitu}
Q(\theta^{n+1}, \eta^{n+1})=Diag(2) \eta^{n+1} - k \varphi(\theta^{n+1}, \eta^{n+1}) - LDQ(\theta^{n}, \eta^{n}),
\end{equation}
where $LDQ = Diag(2) \eta^{n} + k \varphi(\theta^{n}, \eta^{n})$ is known at every instant of time.

Furthermore:
\begin{equation} \label{diag2}
Diag(2)=\left[
\begin{array}{cccccccc}
2 & 0  &    0      &   0      &  \cdots  &     0    &    0     &    0   \\
0  & 2 & 0  &   0      &  \cdots  &     0    &    0     &    0   \\
   0      &  0 & 2 & 0 &          &     0    &    0     &    0   \\
 \vdots   &           &           &          &  \ddots  &          &          &  \vdots \\
     0    &    0      &    0      &    0     &          & 0 &2 & 0 \\
     0    &    0      &    0      &    0     &          &   0      & 0 & 2      
     \end{array}
\right],
\end{equation}
\vspace{0.8cm}
\begin{equation}\label{fi}
\varphi^{n}=\varphi(\theta ,\eta^{n})=\left[
 \begin{array}{ccc}
\varphi^{n}_{1} \\
\varphi^{n}_{2} \\
\varphi^{n}_{3} \\
  \vdots \\
\varphi^{n}_{M-1}  \\ 
\varphi^{n}_{M} \\
     \end{array}
\right], \  \ \ 
\eta^{n}=\left[
 \begin{array}{ccc}
\eta^{n}_{1} \\
\eta^{n}_{2} \\
\eta^{n}_{3} \\
  \vdots \\
\eta^{n}_{M-1}  \\ 
\eta^{n}_{M} \\
     \end{array}
\right].
\end{equation}
Therefore, the discrete form of (\ref{F en situ}) and (\ref{Q en situ}) is given by (\ref{Fdiscre ensitu}) and (\ref{Qdiscre ensitu}).
\begin{eqnarray} \label{FQ ensitu}
G^{n}(\theta^{n+1}, \eta^{n+1}) \bullet \theta^{n+1} &=& 0, \\
Q^{n}(\theta^{n+1}, \eta^{n+1}) &=& 0,
\end{eqnarray}

and must comply with:

\begin{equation} \label{theta ensitu}
\theta^{n+1} \geq 0.
\end{equation}

Thus, joining  (\ref{Fdiscre ensitu}), (\ref{Qdiscre ensitu}) (\ref{FQ ensitu}) and (\ref{theta ensitu}) form a Mixed Complementarity Problem, which can be solved by the FDA-MNCP algorithm

\begin{algorit}[implementation FDA-MNCP.]\label{implementaofda-mncp para insitu} 

\item[\textbf{Step 1.}] Do $n=0$ and $N=1 / \Delta t.$  

\item[\textbf{Step 2.}] To obtain $\theta^{n+1}$ and $\eta^{n+1}$ we apply the method $FDA-MNCP$ to solve the mixed complementarity problem.
\begin{eqnarray}
G^{n}(\theta^{n+1}, \eta^{n+1}) \bullet \theta^{n+1} = 0, \ \ \theta^{n+1} \geq 0, \nonumber \\
Q^{n}(\theta^{n+1}, \eta^{n+1}) = 0
\end{eqnarray}
with
\begin{eqnarray} 
G^{n}(\theta^{n+1}, \eta^{n+1}) = A \theta^{n+1} + \lambda P(\theta^{n+1}, \eta^{n+1}) - 2k \Phi(\theta^{n+1}, \eta^{n+1})-LD(\theta^{n}, \eta^{n+1}) \geq 0, \nonumber \\
Q(\theta^{n+1}, \eta^{n+1})=Diag(2)\cdot \eta^{n+1} - k \varphi(\theta^{n+1}, \eta^{n+1}) - LDQ(\theta^{n}, \eta^{n+1}) = 0. \nonumber
\end{eqnarray}
where the matrices $A, B$ and the vectors $P, \Phi$ and $\varphi$ are given in (\ref{A ensitu}), (\ref{B ensitu}), (\ref{P e phi ensitu}), and (\ref{fi}). 

\item[\textbf{Step 3.}] If $n = N$ then \textbf{END.} \\
 \ \  \ \ \ \ \  \ \ else $n = n + 1$ return to Step \textbf{2}

\end{algorit}

The numerical results obtained from the implementation of the Algorithm in Matlab are shown in the following Section.

\section{Comparison of FDA-MNCP and FDA-NCP methods}
We will now present the numerical results of the simulations performed in Matlab
for the FDA-MNCP method. For this simulation, we considered $[x_{0}, x_{M}]=[0, \ 0.05]$ and $[t_{0}, t_{N}]=[0, 1]$ as the space and time intervals respectively. We keep the numbers of subintervals in time constant and equal to: $N=10^5$, that is,  $k= \Delta t = 10^{-5}$ while the number of subintervals in space will be equal to $M=50, 100, 200, 400$. For
the FDA-MNCP method we consider an error tolerance of $10^{-8}$. 

\begin{figure}[H] 
\centering
\includegraphics[scale=0.65]{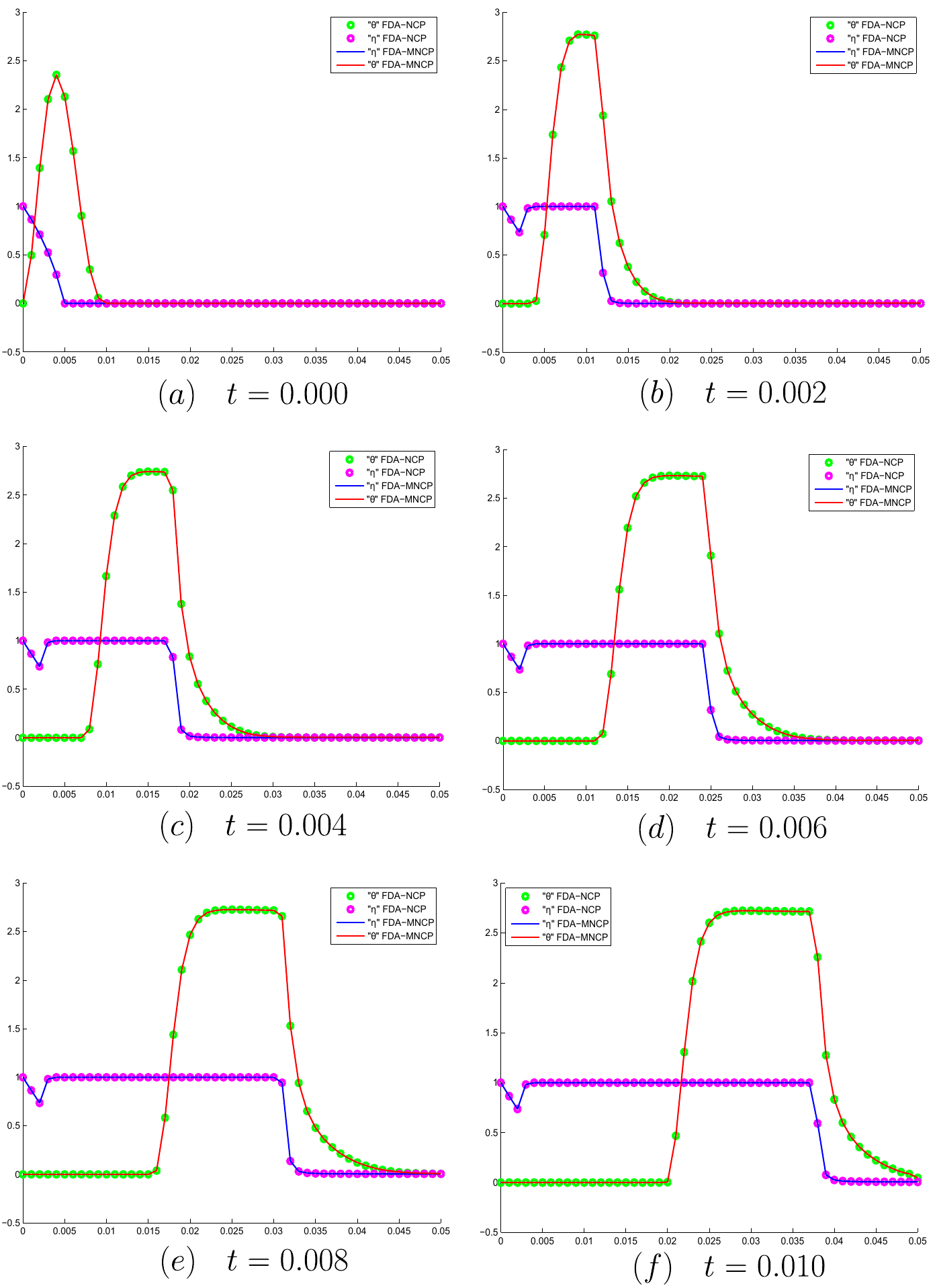} 
\caption{Comparison of the FDA-MNCP and FDA-NCP methods for M = 50 at time instants $t = 0.000, 0.002, 0.004, 0.006, 0.008, 0.010$. The values of $\theta$ are represented by green dots and a solid red line, the values of $\eta$ are represented by pink dots and a solid blue line.} \label{MNCP1-H50}
\end{figure}

\begin{figure}[H] 
\centering
\includegraphics[scale=0.65]{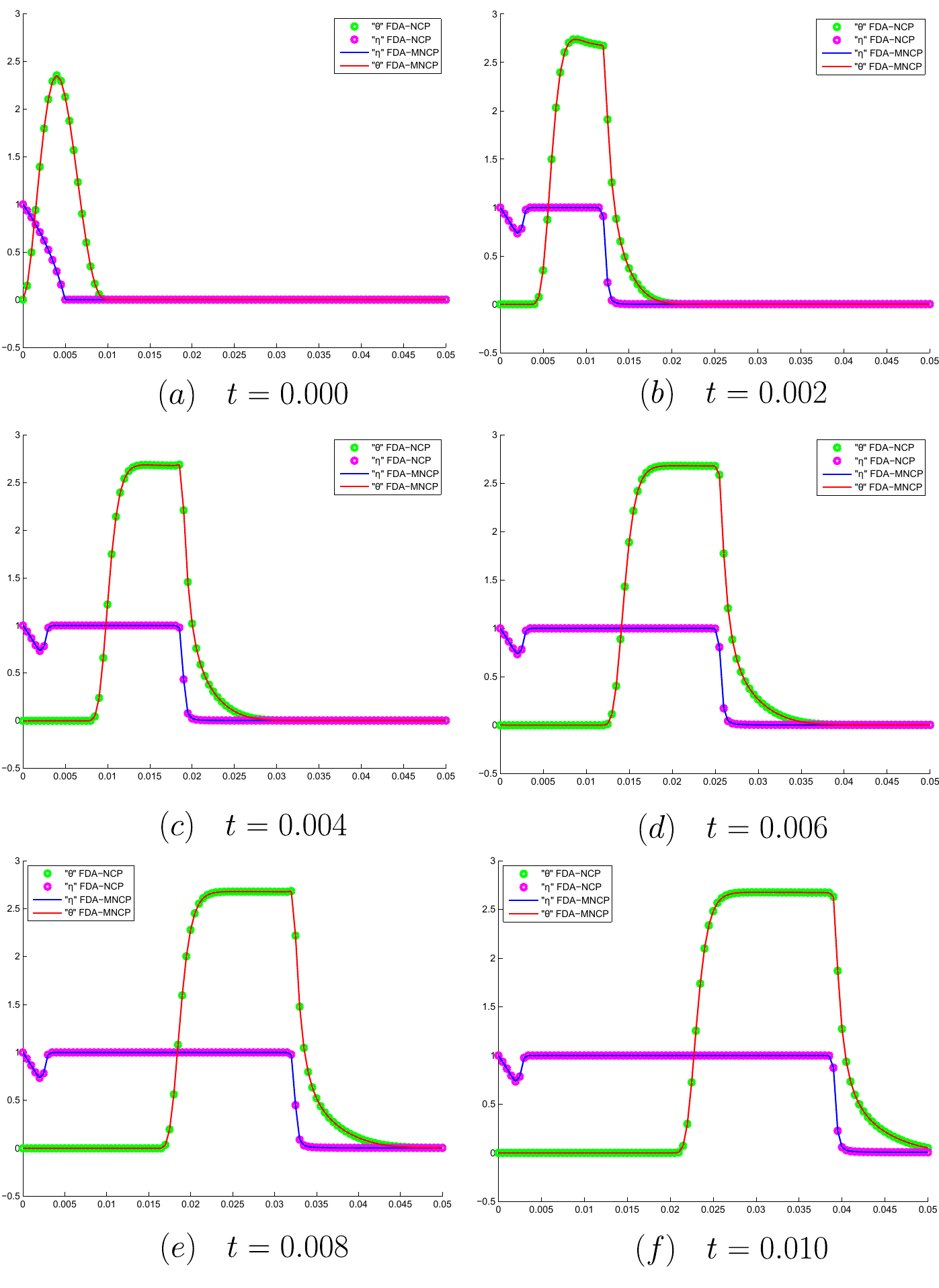}
\caption{Comparison of the FDA-MNCP and FDA-NCP methods for$ M = 100$ at times $t = 0.000, 0.002, 0.004, 0.006, 0.008, 0.010$. The values of $\theta$ are represented by green dots and a solid red line, the values of $\eta$ are represented by pink dots and a solid blue line.}  \label{MNCP1-H100}

\end{figure}

\begin{figure}[H]
\centering
\includegraphics[scale=0.65]{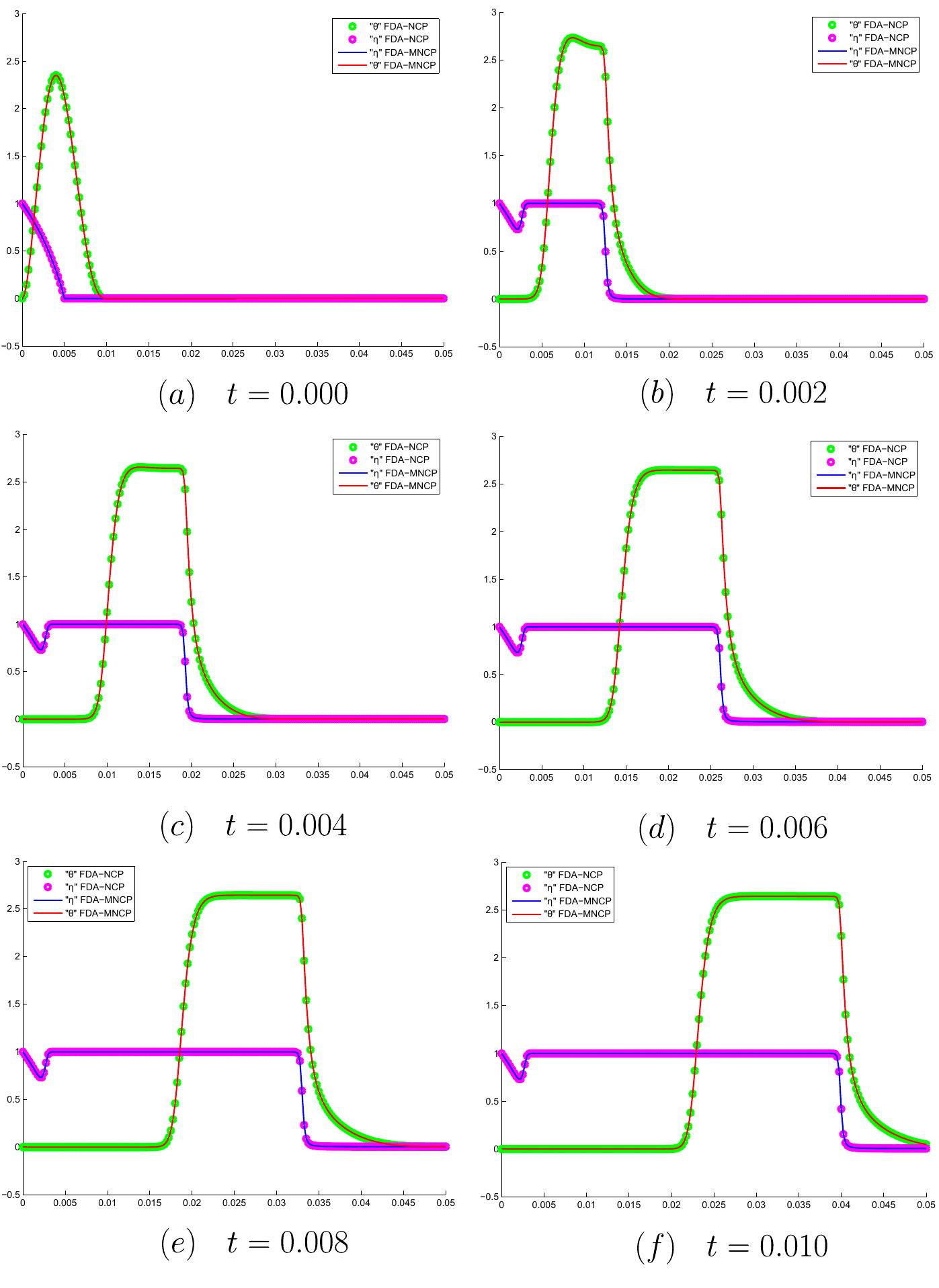} 
\caption{Comparison of the FDA-MNCP and FDA-NCP methods for $M = 200$ at times $t = 0.000, 0.002, 0.004, 0.006, 0.008, 0.010$. The values of $\theta$ are represented by green dots and a solid red line, the values of $\eta$ are represented by pink dots and a solid blue line.}  \label{MNCP1-H200}

\end{figure}

\begin{figure}[H] 
\centering
\includegraphics[scale=0.65]{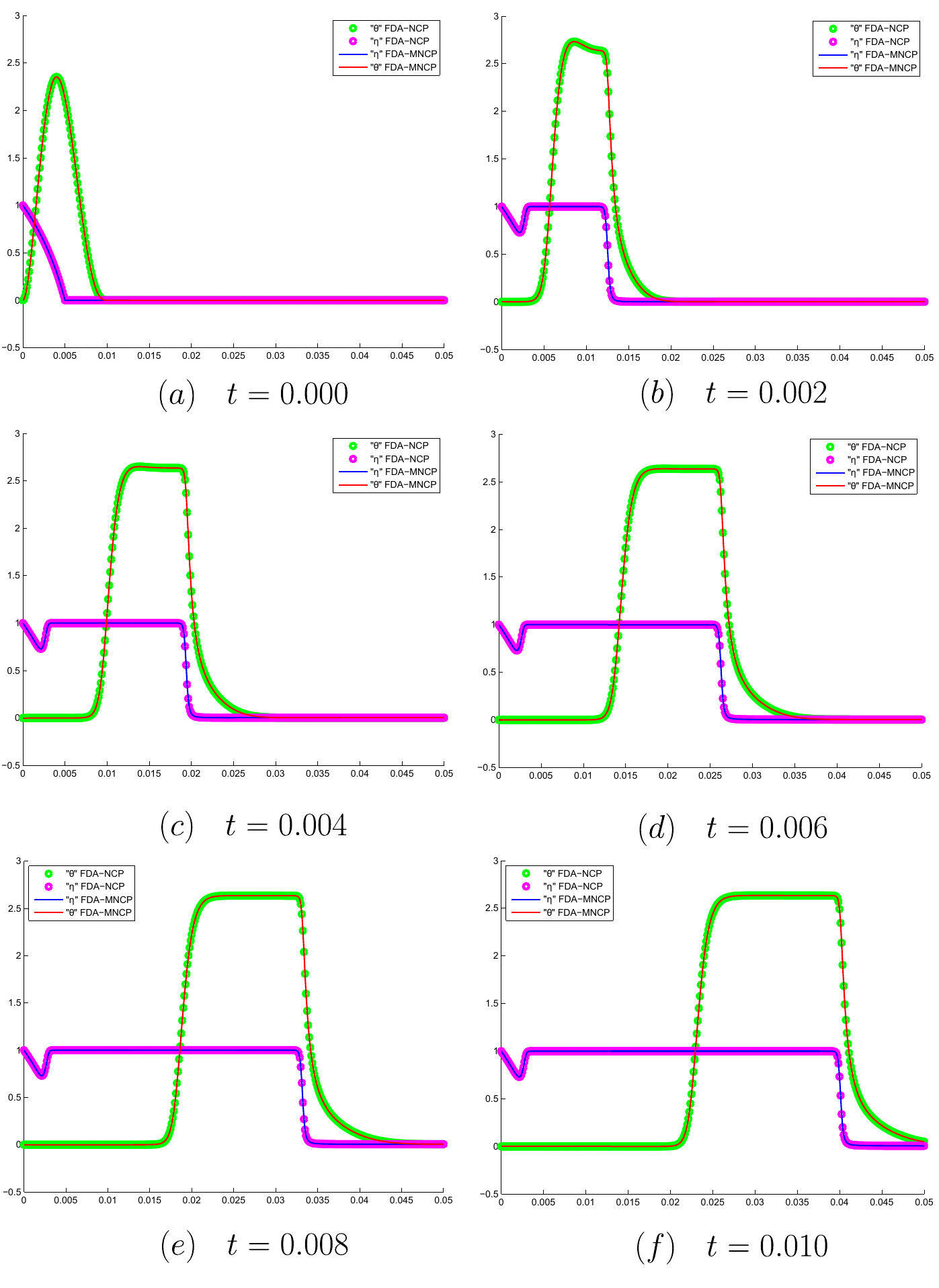} 
\caption{Comparison of the FDA-MNCP and FDA-NCP methods for $M = 400$ at times $t = 0.000, 0.002, 0.004, 0.006, 0.008, 0.010$. The values of $\theta$ are represented by green dots and a solid red line, the values of $\eta$ are represented by pink dots and a solid blue line. } \label{MNCP1-H400}
\end{figure}

The values of the dimensionless parameters in (\ref{constantesadimensionais}) are: 
\begin{eqnarray}
x^{\star} = 9,1 \times 10^{4}  [m], \ \ \ \ t^{\star} = 1,48 \times 10^{8}  [s], \ \ \ \Delta T^{\star} = 74,4 [K], \nonumber \\   u^{\star} = 6,1 \times 10^{-4}, \ \ \ \ \ \ \ \ P_{e_{T}} = 1406, \ \ \ \ \ \ \ \ \ \beta = 7,44 \times 10^{10}, \nonumber \\
E = 93,8 \ \ \ \ \ \ \ \ \ \ \ \ \ \ \ \ \ \ \ \ \  \theta_{0} = 3.67 \ \ \ \ \ \ \ \ \ \ \ \ \ \ \ \ \ \ u=3,76 . \nonumber 
\end{eqnarray}
with the previous input data we obtain Figures (\ref{MNCP1-H50}),  (\ref{MNCP1-H100}), (\ref{MNCP1-H200}), (\ref{MNCP1-H400}), which show the results obtained by algorithm \ref{implementaofda-mncp para insitu} and algorithm 5 of  \cite{angel} for the FDA-MNCP and FDA-NCP  \cite{sandrofdancp} methods, respectively.

\newpage

\newpage

\newpage

\begin{figure}[H] 
\centering
\includegraphics[scale=0.4]{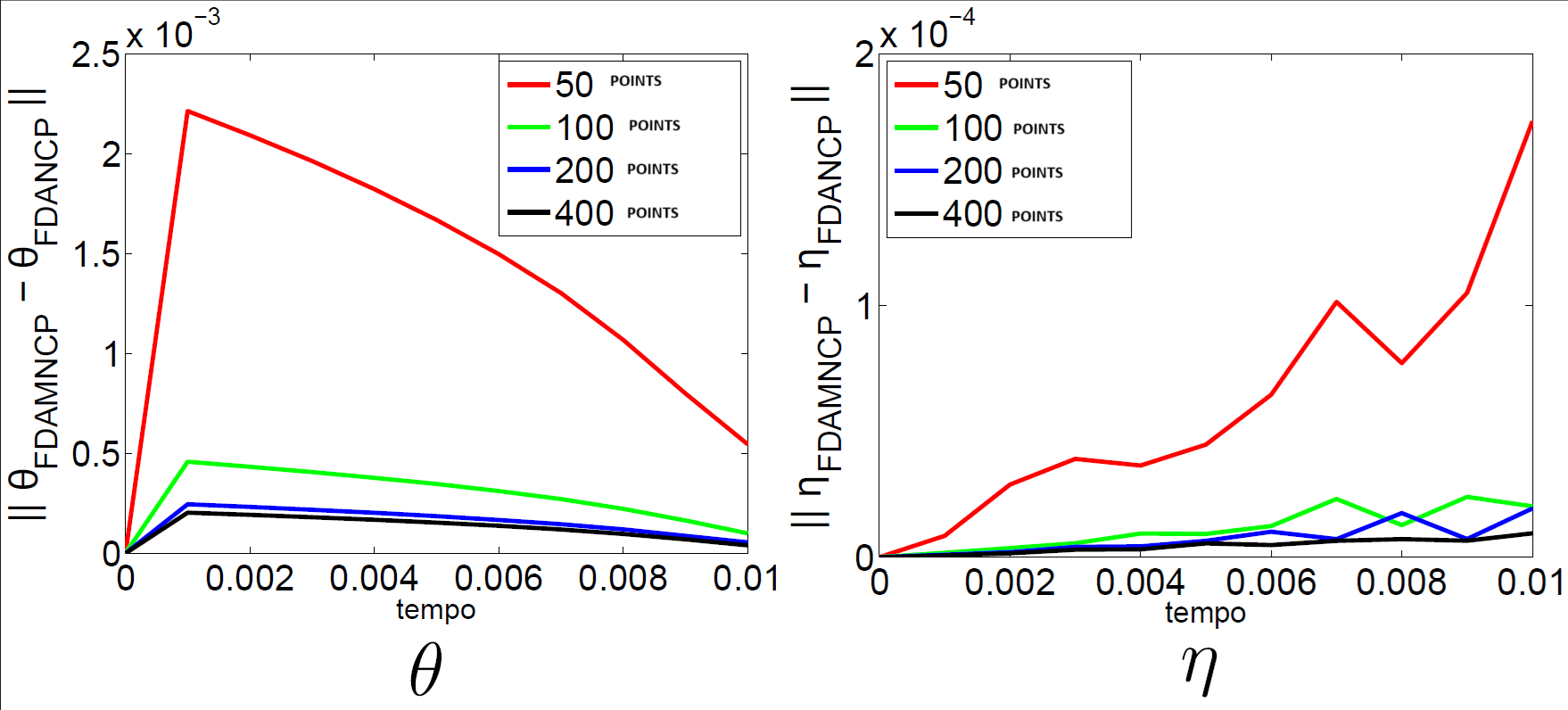} 
\caption{Difference between DA-NCP and FDA-NCP methods.} \label{errorMNCP-NCP theta eta}
\end{figure}

As we can see in Figures \ref{MNCP1-H50}, \ref{MNCP1-H100}, \ref{MNCP1-H200} and \ref{MNCP1-H400} obtained previously, they show us that the results obtained by the FDA-MNCP methods and FDA-NCP, coincide very well, as can be seen in the times indicated in the figures presented previously. In Figure \ref{errorMNCP-NCP theta eta}, the differences between the methods can be observed.
we can see that the difference between the solutions of $\theta$ and $\eta$ are very small
as we increase the number of points, between the FDA-MNCP and FDA-NCP methods\cite{sandrofdancp}.

Below we show four tables where we compare the computational process time for  $M=50, \ 100, \ 20, \ 400$ of the FDA-MNCP method and the FDA-NCP method studied in \cite{angel}.
 \begin{table}[H]
\begin{center} 
\begin{tabular}{|c|ccccc|ccccc|} \hline
 \ \  & \multicolumn{5}{c|}{FDA-MNCP} & \multicolumn{5}{c|}{FDA-NCP} \\ 
  $t$&$t(n)$&iter&BL(t)&$[S]$&$[\nabla S]$&$t(n)$&iter&BL(t)&$F$&$\nabla F$  \\ \hline
0.001 & 0.401 & 33 & 1.0 & 43 & 34 & 0.194 & 20 & 1.0 & 43 & 40 \\
0.002 & 0.411 & 34 & 0.8 & 44 & 35 & 0.194 & 20 & 1.0 & 43 & 40 \\
0.003 & 0.406 & 34 & 1.0 & 43 & 35 & 0.199 & 20 & 1.0 & 43 & 40 \\
0.004 & 0.386 & 32 & 0.8 & 41 & 33 & 0.184 & 20 & 1.0 & 43 & 40 \\
0.005 & 0.385 & 32 & 0.8 & 41 & 33 & 0.178 & 21 & 1.0 & 45 & 42 \\
0.006 & 0.406 & 34 & 0.8 & 43 & 35 & 0.198 & 21 & 1.0 & 45 & 42 \\
0.007 & 0.395 & 33 & 0.8 & 42 & 34 & 0.193 & 21 & 1.0 & 45 & 42 \\
0.008 & 0.403 & 33 & 0.8 & 43 & 34 & 0.229 & 21 & 1.0 & 45 & 42 \\
0.009 & 0.390 & 32 & 0.8 & 42 & 33 & 0.169 & 21 & 1.0 & 45 & 42 \\
0.010 & 0.405 & 33 & 0.8 & 44 & 34 & 0.173 & 21 & 1.0 & 45 & 42 \\ \hline
\end{tabular}
\caption{Comparison of the computational process time with $M = 50$ for the FDA-MNCP and FDA-NCP methods. $t(n)$ is the time measured in seconds that the method took to find the solution at time $t$.}  \label{tabela de h=50}
\end{center}
\end{table}
As we can see in Table  \ref{tabela de h=50}, the computational process time used by the FDA-MNCP method doubles the time of the FDA-NCP method for a partition of the $x$-axis into fifty points; the number of iterations of the FDA-MNCP method for the indicated times is also greater than that of the FDA-NCP method. In addition, the third column of each table shows the value of $t^k$ for each of the methods.

\begin{table}[H]
\begin{center}
\begin{tabular}{|c|ccccc|ccccc|} \hline
 \ \  & \multicolumn{5}{c|}{FDA-MNCP} & \multicolumn{5}{c|}{FDA-NCP} \\ 
  $t$&$t(n)$&iter&BL(t)&$[S]$&$[\nabla S]$&$t(n)$&iter&BL(t)&$F$&$\nabla F$  \\ \hline
0.001 & 0.912 & 36 & 0.8  & 47 & 37 & 0.555 & 21 & 1.0 & 45 & 42 \\
0.002 & 0.885 & 36 & 0.8  & 46 & 37 & 0.519 & 21 & 1.0 & 45 & 42 \\
0.003 & 0.870 & 37 & 0.8  & 47 & 38 & 0.515 & 21 & 1.0 & 45 & 42 \\
0.004 & 0.892 & 35 & 0.64 & 46 & 36 & 0.494 & 21 & 1.0 & 45 & 42 \\
0.005 & 0.831 & 34 & 0.8  & 44 & 35 & 0.486 & 21 & 1.0 & 45 & 42 \\
0.006 & 0.793 & 34 & 0.8  & 43 & 35 & 0.508 & 21 & 1.0 & 45 & 42 \\
0.007 & 0.790 & 34 & 0.8  & 43 & 35 & 0.457 & 21 & 1.0 & 45 & 42 \\
0.008 & 0.894 & 36 & 0.64 & 49 & 37 & 0.515 & 21 & 1.0 & 45 & 42 \\
0.009 & 0.869 & 35 & 0.8  & 46 & 36 & 0.507 & 21 & 1.0 & 45 & 42 \\
0.010 & 0.839 & 35 & 0.8  & 46 & 36 & 0.479 & 21 & 1.0 & 45 & 42 \\ \hline
\end{tabular}
\caption{Comparison of the computational process time with $M = 100$ for the FDA-MNCP and FDA-NCP methods. t(n) is the time measured in seconds that the method took to find the solution at time $t$.}  \label{tabela de h=100} 
\end{center}
\end{table}
In Table \ref{tabela de h=100}, we can also observe that the computational process time
used by the FDA-MNCP method is slightly greater than the time of the FDA-NCP method
as the number of iterations.

 \begin{table}[H]
\begin{center}
\begin{tabular}{|c|ccccc|ccccc|} \hline
 \ \  & \multicolumn{5}{c|}{FDA-MNCP} & \multicolumn{5}{c|}{FDA-NCP} \\ 
  $t$&$t(n)$&iter&BL(t)&$[S]$&$[\nabla S]$&$t(n)$&iter&BL(t)&$F$&$\nabla F$  \\ \hline
0.001 & 1.771 & 36 & 0.8 & 47 & 37 & 3.223 & 21 & 1.0 & 45 & 42 \\
0.002 & 1.852 & 38 & 0.8 & 48 & 39 & 3.185 & 21 & 1.0 & 45 & 42 \\
0.003 & 1.969 & 37 & 0.8 & 46 & 38 & 3.129 & 21 & 1.0 & 45 & 42 \\
0.004 & 1.995 & 38 & 0.8 & 47 & 39 & 3.245 & 21 & 1.0 & 45 & 42 \\
0.005 & 1.771 & 37 & 0.8 & 46 & 38 & 3.143 & 21 & 1.0 & 45 & 42 \\
0.006 & 1.785 & 37 & 1.0 & 45 & 38 & 3.122 & 21 & 1.0 & 45 & 42 \\
0.007 & 1.787 & 37 & 0.8 & 47 & 38 & 3.205 & 21 & 1.0 & 45 & 42 \\
0.008 & 1.778 & 37 & 0.8 & 47 & 38 & 3.241 & 21 & 1.0 & 45 & 42 \\
0.009 & 2.001 & 36 & 0.8 & 46 & 37 & 3.346 & 22 & 1.0 & 47 & 44 \\
0.010 & 1.685 & 35 & 0.8 & 45 & 36 & 3.396 & 22 & 1.0 & 47 & 44 \\ \hline
\end{tabular}
\caption{Comparison of the computational process time with $M = 200$ for the FDA-MNCP and FDA-NCP methods. t(n) is the time measured in seconds that the method took to find the solution at time $t$.}    \label{tabela de h=200}
\end{center}
\end{table}

In Table  \ref{tabela de h=200}, we see that when we increase the number of partitions to $200$, the computational process time of the FDA-MNCP method starts to be shorter than the time of the 1FDA-NCP method, although the number of iterations is still greater.
\begin{table}[H]
\begin{center}
\begin{tabular}{|c|ccccc|ccccc|} \hline
 \ \  & \multicolumn{5}{c|}{FDA-MNCP} & \multicolumn{5}{c|}{FDA-NCP} \\ 
  $t$&$t(n)$&iter&BL(t)&$[S]$&$[\nabla S]$&$t(n)$&iter&BL(t)&$F$&$\nabla F$  \\ \hline
0.001 & 4.497 & 37 & 0.8  & 48 & 38 & 33.127 & 21 & 1.0 & 45 & 21 \\
0.002 & 4.337 & 37 & 1.0  & 46 & 38 & 33.174 & 21 & 1.0 & 45 & 21 \\
0.003 & 4.533 & 38 & 0.8  & 47 & 39 & 33.034 & 21 & 1.0 & 45 & 21 \\
0.004 & 4.494 & 39 & 0.8  & 48 & 40 & 33.398 & 21 & 1.0 & 45 & 21 \\
0.005 & 4.429 & 37 & 1.0  & 45 & 38 & 34.722 & 22 & 1.0 & 47 & 22 \\
0.006 & 4.329 & 37 & 1.0  & 45 & 38 & 34.652 & 22 & 1.0 & 47 & 22 \\
0.007 & 4.523 & 39 & 0.8  & 49 & 40 & 34.777 & 22 & 1.0 & 47 & 22 \\
0.008 & 4.530 & 39 & 0.8  & 49 & 40 & 34.802 & 22 & 1.0 & 47 & 22 \\
0.009 & 4.459 & 37 & 0.8  & 47 & 38 & 34.824 & 22 & 1.0 & 47 & 22 \\
0.010 & 4.315 & 37 & 1.0  & 48 & 38 & 34.571 & 22 & 1.0 & 47 & 22 \\ \hline
\end{tabular}
\caption{Comparison of the computational process time with $M = 400$ for the FDA-MNCP and FDA-NCP methods. $t(n)$ is the time measured in seconds that the method took to find the solution at time $t$.}    \label{tabela de h=400}
\end{center}
\end{table}

With the last Table \ref{tabela de h=400}, we can observe that by considerably increasing the number of partitions, the computational process time of the FDA-MNCP method is much lower than that of the FDA-NCP\cite{sandrofdancp} method, because for the calculation of $\nabla S$ in the FDA-MNCP method, it makes half the calculations than for $\nabla F$ in the FDA-NCP\cite{sandrofdancp} method, although the number of iterations is greater. We can also observe that the number of iterations varies little in each of the tables presented previously.

\section{Error analysis}
\begin{figure}[H]
\centering
\includegraphics[scale=0.6]{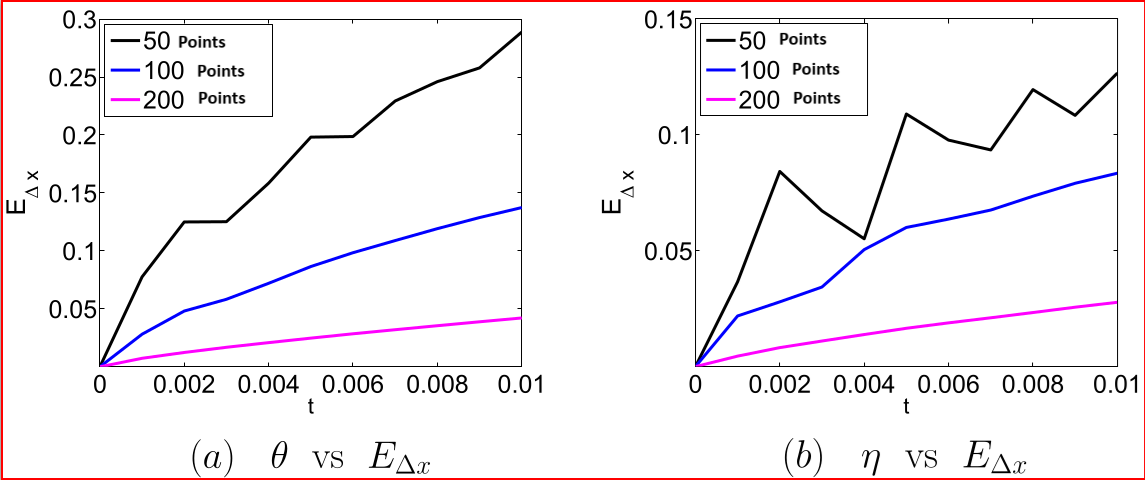} 
\caption{Time $(t) \mbox{vs} E_{\Delta x}$. Relative Method Error. Here $\Delta x = \frac{1}{50}, \ \frac{1}{100}, \ \frac{1}{200}$.} \label{errorelativoFDAMNCP}
\end{figure}

In this section we will perform a numerical study of the relative error for the FDA-MNCP method, for each instant of time and $E_{\Delta x}$, where $\Delta x \ = \ h, \ \ h/2, \ \ h/4$. The length of the subinterval in time will be constant and equal to $\Delta t \ = \ k \ = \ 10^{-5}$ \ e \ $h \ = \ \frac{1}{50}$, and compared with the FDA-NCP method, \cite{angel}. Tables  \ref{tabelaFDAMNCPH=1sobre50theta} and \ref{tabelaFDAMNCPH=1sobre50eta}  show the results for the Relative Errors of the FDA-MNCP method, while Tables \ref{FDANCPH=1sobre50theta} and \ref{FDANCPH=1sobre50eta} show the Relative Error for the FDA-NCP method and are represented in Figures \ref{errorelativoFDAMNCP} and \ref{errorelativoFDANCP}.

\begin{figure}[H]
\centering
\label{errorelativoFDANCP}
\includegraphics[scale=0.37]{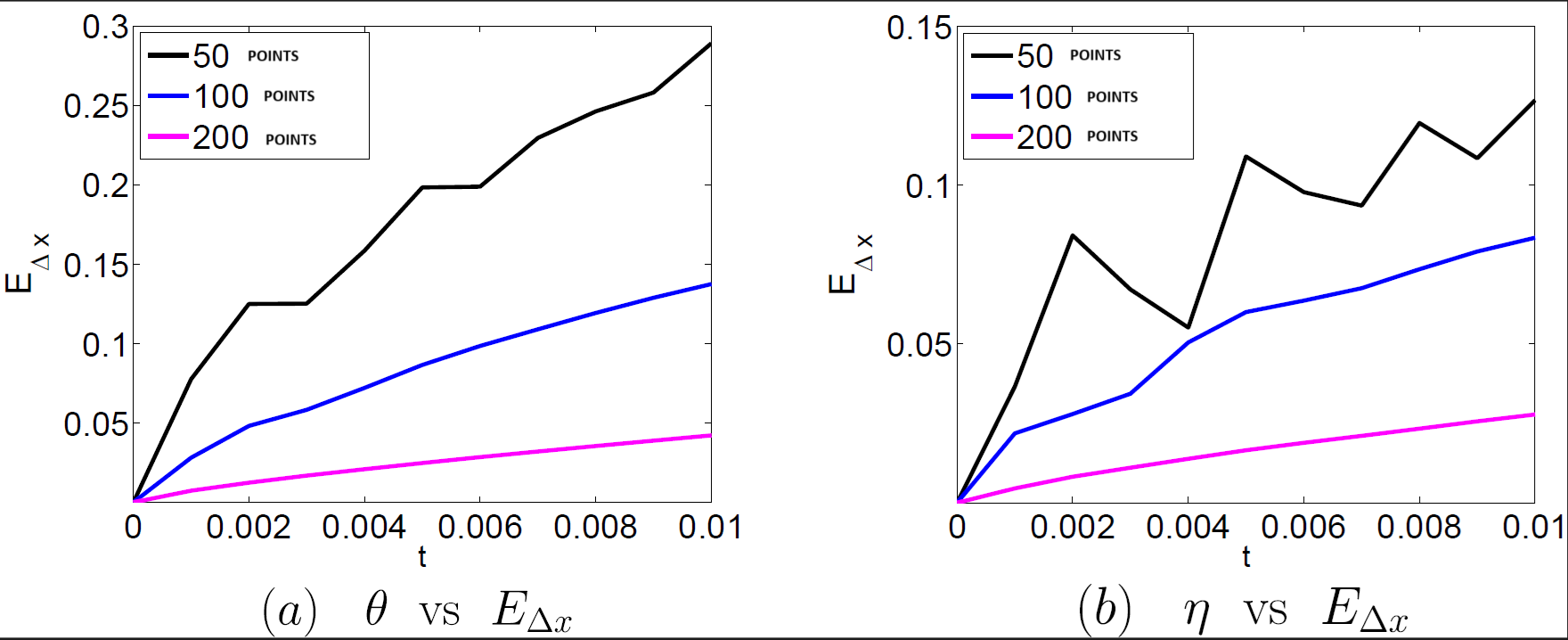} 
\caption{Time $(t) \mbox{vs} E_{\Delta x}$. Relative Method Error FDA-NCP. Here $\Delta x = \frac{1}{50}, \ \frac{1}{100}, \ \frac{1}{200}$.} 
\label{errorelativoFDANCP}
\end{figure}

As can be seen in Figures \ref{errorelativoFDAMNCP} and \ref{errorelativoFDANCP}, the relative errors corresponding to the FDA-MNCP and FDA-NCP methods, respectively, are very similar, as can be seen in the following tables, which show the relative errors for $\theta$ and $\eta$ usando $h = 1/50$.
\begin{table}[H] 
\begin{center}  
\begin{tabular}{c|c|c|c|c|c} \hline
 $t$ & $E_{h}$ & $E_{\frac{h}{2}}$ & $E_{\frac{h}{4}}$ & $\frac{E_{h}}{E_{\frac{h}{2}}}$ & $\frac{E_{\frac{h}{2}}}{E_{\frac{h}{4}}}$ \\ \hline
  0.00100000 & 0.07757130 & 0.02809022 & 0.00730451 & 2.76 & 3.85 \\
  0.00200000 & 0.12495960 & 0.04805076 & 0.01233415 & 2.60 & 3.90 \\
  0.00300000 & 0.12512590 & 0.05824369 & 0.01679689 & 2.15 & 3.47 \\
  0.00400000 & 0.15843617 & 0.07200921 & 0.02080162 & 2.20 & 3.46 \\
  0.00500000 & 0.19817638 & 0.08645161 & 0.02467447 & 2.29 & 3.50 \\
  0.00600000 & 0.19866653 & 0.09841438 & 0.02840411 & 2.02 & 3.46 \\
  0.00700000 & 0.22936049 & 0.10892940 & 0.03195340 & 2.11 & 3.41 \\
  0.00800000 & 0.24606770 & 0.11916519 & 0.03538988 & 2.06 & 3.37 \\
  0.00900000 & 0.25789027 & 0.12872004 & 0.03877097 & 2.00 & 3.32 \\
  0.01000000 & 0.28886846 & 0.13735367 & 0.04209700 & 2.10 & 3.26 \\ \hline
\end{tabular}
\caption{Relative error for $\theta$ com FDA-MNCP e $h= \frac{1}{50}$ for the time instants $t$ indicated in the first column.} \label{tabelaFDAMNCPH=1sobre50theta}
\end{center}
\end{table}
\begin{table}[H]
\begin{center}  
\begin{tabular}{c|c|c|c|c|c} \hline
 $t$ & $E_{h}$ & $E_{\frac{h}{2}}$ & $E_{\frac{h}{4}}$ & $\frac{E_{h}}{E_{\frac{h}{2}}}$ & $\frac{E_{\frac{h}{2}}}{E_{\frac{h}{4}}}$ \\ \hline
0.001 & 0.07753953 & 0.02808885 & 0.00730328 & 2.76 & 3.85 \\
0.002 & 0.12490382 & 0.04804832 & 0.01233527 & 2.60 & 3.89 \\
0.003 & 0.12507648 & 0.05824052 & 0.01679940 & 2.15 & 3.47 \\
0.004 & 0.15841526 & 0.07200834 & 0.02080540 & 2.20 & 3.46 \\
0.005 & 0.19818512 & 0.08645353 & 0.02467944 & 2.29 & 3.50 \\
0.006 & 0.19869586 & 0.09841770 & 0.02841001 & 2.02 & 3.46 \\
0.007 & 0.22940100 & 0.10893437 & 0.03196011 & 2.10 & 3.41 \\
0.008 & 0.24613384 & 0.11917200 & 0.03539747 & 2.07 & 3.37 \\
0.009 & 0.25797041 & 0.12872819 & 0.03877948 & 2.00 & 3.32 \\
0.010 & 0.28897837 & 0.13736335 & 0.04210642 & 2.10 & 3.26 \\ \hline
\end{tabular}
\caption{Relative error for $\theta$ com FDA-NCP e $h= \frac{1}{50}$ for the time instants $t$ indicated in the first column.} \label{FDANCPH=1sobre50theta}
\end{center}
\end{table}

In Tables \ref{tabelaFDAMNCPH=1sobre50theta} and \ref{FDANCPH=1sobre50theta} we see that the relative errors $E_{h}, \ E_{\frac{h}{2}}, \ E_{\frac{h}{4}} \ $  for  $\ \theta \ $ of both the FDA-MNCP method and the FDA-NCP method are very similar since they are the same in the first three decimal places. 

\begin{table}[H] 
\begin{center}  
\begin{tabular}{c|c|c|c|c|c} \hline
 $t$ & $E_{h}$ & $E_{\frac{h}{2}}$ & $E_{\frac{h}{4}}$ & $\frac{E_{h}}{E_{\frac{h}{2}}}$ & $\frac{E_{\frac{h}{2}}}{E_{\frac{h}{4}}}$ \\ \hline
  0.00100000 & 0.03656879 & 0.02183277 & 0.00451877 & 1.67 & 4.83 \\
  0.00200000 & 0.08417410 & 0.02792488 & 0.00819559 & 3.01 & 3.41 \\
  0.00300000 & 0.06717364 & 0.03432880 & 0.01103677 & 1.96 & 3.11 \\
  0.00400000 & 0.05513563 & 0.05046208 & 0.01381861 & 1.09 & 3.65 \\
  0.00500000 & 0.10897196 & 0.06002542 & 0.01653614 & 1.82 & 3.63 \\
  0.00600000 & 0.09771893 & 0.06361601 & 0.01887177 & 1.54 & 3.37 \\
  0.00700000 & 0.09347666 & 0.06753455 & 0.02105732 & 1.38 & 3.21 \\
  0.00800000 & 0.11948456 & 0.07350345 & 0.02332175 & 1.63 & 3.15 \\
  0.00900000 & 0.10838295 & 0.07908666 & 0.02560247 & 1.37 & 3.09 \\
  0.01000000 & 0.12663557 & 0.08343716 & 0.02775604 & 1.52 & 3.01 \\ \hline
\end{tabular}
\caption{Relative error for $\eta$ with FDA-MNCP and $h = \frac{1}{50}$ for the instants of time $t$ indicated in the first column.} \label{tabelaFDAMNCPH=1sobre50eta}
\end{center}
\end{table}

\begin{table}[H]
\begin{center}  
\begin{tabular}{c|c|c|c|c|c} \hline
 $t$ & $E_{h}$ & $E_{\frac{h}{2}}$ & $E_{\frac{h}{4}}$ & $\frac{E_{h}}{E_{\frac{h}{2}}}$ & $\frac{E_{\frac{h}{2}}}{E_{\frac{h}{4}}}$ \\ \hline
0.001 & 0.03656752 & 0.02183277 & 0.00451876 & 1.67 & 4.83 \\
0.002 & 0.08415103 & 0.02792385 & 0.00819578 & 3.01 & 3.41 \\
0.003 & 0.06714517 & 0.03432836 & 0.01103751 & 1.96 & 3.11 \\
0.004 & 0.05513010 & 0.05046448 & 0.01382034 & 1.09 & 3.65 \\
0.005 & 0.10898570 & 0.06002658 & 0.01653867 & 1.82 & 3.63 \\
0.006 & 0.09776211 & 0.06361489 & 0.01887480 & 1.54 & 3.37 \\
0.007 & 0.09351854 & 0.06753575 & 0.02106120 & 1.38 & 3.21 \\
0.008 & 0.11952653 & 0.07350797 & 0.02332674 & 1.63 & 3.15 \\
0.009 & 0.10846255 & 0.07909222 & 0.02560830 & 1.37 & 3.09 \\
0.010 & 0.12671978 & 0.08344353 & 0.02776237 & 1.52 & 3.01 \\ \hline
\end{tabular}
\caption{Relative error for $\eta$ with FDA-NCP and $h= \frac{1}{50}$ for the instants of time $t$ indicated in the first column.} \label{FDANCPH=1sobre50eta}
\end{center}
\end{table}

Similarly, in Tables \ref{tabelaFDAMNCPH=1sobre50eta} and \ref{FDANCPH=1sobre50eta} we $E_{h}, \ E_{\frac{h}{2}}, \ E_{\frac{h}{4}} \ $ for $\eta$ of both the FDA-MNCP method and the FDA-NCP method are very similar,
since the first three decimal places are the same
\section{Conclusions}
\begin{itemize}
\item We propose the solution of the simple in-situ combustion model using a numerical method based on an implicit finite difference scheme and a nonlinear mixed complementarity algorithm, which can be applied to parabolic problems and can be rewritten in the form of a mixed complementarity problem. This method is applied to the system(\ref{S1}). 

\item In this work, we have seen that the feasible interior point algorithm FDA-MNCP, as a good technique to numerically solve mixed complementarity problems. Theoretical results of the algorithm ensure the global convergence. 

\item We solve the in-situ combustion model(\ref{S1}) using the nonlinear mixed complementarity algorithm and comparing it to the FDA-NCP method, with the two solutions being very close, as can be seen in Figures  \ref{MNCP1-H50}, \ref{MNCP1-H100}, \ref{MNCP1-H200} and \ref{MNCP1-H400}. This suggests that it is possible to apply this method to parabolic and hyperbolic problems, which can be written as a mixed complementarity problem. The FDA-MNCP method shows the advantage of being able to be used with more points in the discretization of the space, as can be seen in Tables \ref{tabela de h=200} and \ref{tabela de h=400}, in addition to being faster in computational time than the FDA-NCP method when we increase the discretization in the space. 
\item Regarding the relative errors, Tables \ref{tabelaFDAMNCPH=1sobre50theta} and \ref{tabelaFDAMNCPH=1sobre50eta}  show good evidence of convergence of the mixed complementarity algorithm, which is observed in Figure \ref{errorMNCP-NCP theta eta}, which indicates a decrease in linear growth when we refine the mesh.
\end{itemize}

\end{document}